\documentclass[11pt]{article}
\usepackage{url}
\usepackage{graphics}

\title{L1Packv2: A Mathematica package for minimizing an $\ell_1$-penalized functional}
\author{Ignace Loris\\Mathematics Department, Vrije Universiteit Brussel,\\ Pleinlaan 2, 1050 Brussel, Belgium}

\begin{document}

\maketitle

\begin{abstract}
L1Packv2 is a Mathematica package that contains a number of
algorithms that can be used for the minimization of an
$\ell_1$-penalized least squares functional. The algorithms can
handle a mix of penalized and unpenalized variables. Several
instructive examples are given. Also, an implementation that yields
an exact output whenever exact data are given is provided.
\end{abstract}

\section{Introduction}

In physical and applied sciences, one often encounters the situation
where the quantities in which one is ultimately interested cannot be
measured directly. Such a situation appears e.g. in magnetic
resonance imaging (Fourier components of an image are measured),
optics (convolution of an object with a non-trivial point spread
function forms an image) etc.

Often, it happens that a linear relationship exists between the data
and the quantities of interest. Nonetheless, several issues stand in
the way of `solving' such \emph{linear inverse problems}.
Frequently, the data are incomplete (underdetermined system from a
mathematical point of view) and contaminated with noise
(inconsistent equations from a mathematical point of view).
Moreover, the linear operator linking both is often singular.

A much-used strategy to arrive at a well-defined solution for the
above mentioned linear equations, is to formulate the problem as the
minimization of a quadratic functional consisting of the sum of the
discrepancy $\|Kx-y\|^2$ and a penalization contribution $\|x\|^2$:
\begin{equation}
\arg\min_x \|Kx-y\|^2+\lambda \|x\|^2\label{l2functional}
\end{equation}
in terms of the noisy data $y$ and the linear operator $K$. The
regularization parameter $\lambda$ controls the trade-off between
the data misfit $y-Kx$ and the $\ell_2$-norm of the model
parameters. With a proper choice of $\lambda$, the influence of
noise (present in the data $y$) on the reconstruction of $x$ can now
be kept under control
\cite{Bertero1989,Engl.Hanke.ea1996,Bertero.Boccacci1998}.

The resulting variational equations
\begin{equation}
K^TKx+\lambda x=K^Ty\label{linsys}
\end{equation}
for the minimization problem (\ref{l2functional}) are linear; a
considerable advantage of this technique is that existing efficient
linear solvers (i.e. thoroughly tested algorithms and computer
software) can be used. A disadvantage of using a penalty term of
type $\lambda \|x\|^2$ is that it is quite generic and does not take
into account any a priori information that one may have about the
desired reconstructed object (other than being of finite size).

Recently, scientist have become interested in using
`sparsity-promoting' penalizations for the solution of linear
ill-posed inverse problems. In particular, the following
minimization problem involving an $\ell_1$-penalized least squares
functional has attracted a great deal of attention:
\begin{equation}
\bar x(\lambda)=\arg\min_x \|Kx-y\|^2+2\lambda
\|x\|_1\label{l1functional}
\end{equation}
for a given real matrix $K$, real data $y$ and positive penalization
parameter $\lambda\geq0$. The $\ell_1$-norm is defined as
$\|x\|_1=\sum_i|x_i|$.

The minimization problem (\ref{l1functional}) is commonly referred
to as `lasso' after \cite{Tibshirani1996}, although earlier uses of
this technique do exist \cite{Santosa.Symes1986}. In the context of
signal decomposition, the term `basis pursuit denoising' is often
used \cite{Chen.Donoho.ea1998}.

\emph{Compressed sensing} \cite{Donoh2006,CaRoT2006}, a rapidly
developing research area, is based on the realization that a signal
of length $N$ can often be fully recovered by far fewer than $N$
measurements if it is known in advance that the signal is
\emph{sparse}. The recovery is done by calculating $\arg\min_{Kx=y}
\|x\|_1$, or in case of the presence of noise by calculating the
minimizer (\ref{l1functional}) for which it is known that $\bar
x(\lambda)$ is typically sparse, i.e. a great number of its
components are exactly zero (at least for a large enough penalty
$\lambda$).

Unfortunately, the variational equations that describe the minimizer
(\ref{l1functional}) are nonlinear (see equations (\ref{kkt})
below). The use of a traditional $\ell_2$ penalization leads to
linear equations but a non-sparse minimizer. Hence the question of
finding ways of efficiently recovering the minimizer of functional
(\ref{l1functional}) has attracted a great
deal of interest \cite{Daubechies.Defrise.ea2004,%
Daubechies.Fornasier.ea2007,Koh.Kim.ea2007,%
Figueiredo.Nowak.ea2008,ChCPW2007}.

The Mathematica package L1Packv2 implements a couple of algorithms
for finding the minimizer (\ref{l1functional}). There are five
iterative algorithms which yield an approximation of $\bar
x(\lambda)$ when stopped after a finite number of steps. The
principal algorithm however can calculate the minimizer $\bar
x(\lambda)$ exactly (up to computer round-off) in a finite number of
steps. This algorithm is known as the `homotopy method'
\cite{Osborne.Presnell.ea2000} or as Least Angle Regression (LARS)
\cite{Efron.Hastie.ea2004}.

The algorithms in this package can also treat the minimization of a
slightly more general functional:
\begin{equation}
\min_x \|Kx-y\|^2+2\lambda \sum_i w_i |x_i|
\end{equation}
with positive (or zero) weights $w_i\geq0$.

The procedures are all written with real matrices and signals
in mind. They do not work for complex variables. Extension to
complex data would require optimization over quadratic cones;
the thresholded Landweber algorithm could be easily adapted,
but not the other algorithms. The package was developed and
tested with Mathematica 5.2 \cite{WolframResearch2005}. Other
toolboxes geared towards the recovery of sparse solutions of
linear equations exist for Matlab: SparseLab
\cite{Donoho.Stodden.ea2007}, $\ell_1$-magic
\cite{Candes.Romberg2005}.

\section{Implementation}

In order to describe the solution of the minimization problem
(\ref{l1functional}), it is worthwhile mentioning the following four
equivalent problems:
\begin{enumerate}
\item the minimization of the penalized
least squares functional:
\begin{equation}
\bar x(\lambda)=\arg\min_x \|Kx-y\|^2+2\lambda \|x\|_1\label{xbar}
\end{equation}
\item the constrained minimization problem:
\begin{equation}
\tilde x(R)=\arg\min_{\|x\|_1\leq R} \|Kx -y\|^2\label{xtilde}
\end{equation}
(with an implicit relation between $R$ and $\lambda$, see below)
\item the solution of the variational equations:
\begin{equation}
\begin{array}{lclcl}
(K^T(y-K\bar x))_i&=&\lambda\; \mathrm{sign}(\bar x_i) & \qquad\mathrm{if} & \bar x_i\neq 0 \\
|(K^T(y-K\bar x))_i|&\leq&\lambda & \qquad\mathrm{if} & \bar x_i= 0
\end{array}\label{kkt}
\end{equation}
\item the solution of the fixed-point equation:
\begin{equation}
\bar x=S_\lambda[\bar x+K^T(y-K\bar x)],\label{fixedpoint}
\end{equation}
where $S_\lambda$ is soft-thresholding applied component-wise (see
expression (\ref{soft})).
\end{enumerate}
One has $R=\|\bar x(\lambda)\|_1$ and also:
\begin{equation}
\lambda  =  \|K^T(y-K \bar x(\lambda))\|_\infty
 =  \max_i |(K^T(y-K \bar x(\lambda)))_i|\label{lambdaexpr}
\end{equation}
which follows immediately from relations (\ref{kkt}).

The equivalence of the four formulations
(\ref{xbar}--\ref{fixedpoint}) can easily be seen as follows:
Constrained minimization (\ref{xtilde}) is turned into problem
(\ref{xbar}) by the introduction of a Lagrange multiplier $\lambda$.
The variational equations (\ref{kkt}) are derived from (\ref{xbar})
by putting $x=\bar x+ h e_i$ (where $h$ is a small parameter and
$e_i$ is the $i$th basis vector) and expressing that $\|K\bar
x-y\|^2+2\lambda \|\bar x\|_1\leq \|Kx-y\|^2+2\lambda \|x\|_1$ for
all $h$. The fixed-point equation (\ref{fixedpoint}) is derived from
equation (\ref{kkt}) by adding $x_i$ to both sides of the
equalities, by rearranging terms and using $S_\lambda(u)=u-\lambda\,
\mathrm{sgn}(u)$ for $|u|\geq \lambda$ and $S_\lambda(u)=0$ for
$|u|\leq \lambda$.

It can be shown that $\bar x(\lambda)$ is a piecewise linear
function of $\lambda$ and that $\bar x(\lambda)=0$ for $\lambda\geq
\lambda_\mathrm{max}\equiv \max_i|(K^Ty)_i|$. Hence, it is
sufficient to calculate the nodes of $\bar x(\lambda)$; for generic
values of $\lambda$ in between these nodes, one can use simple
linear interpolation to determine $\bar x(\lambda)$.

The nodes of  $\bar x(\lambda)$ can be calculated using only a
finite number of elementary operations (addition, multiplication,
\ldots) using the LARS/homotopy method
\cite{Osborne.Presnell.ea2000,Efron.Hastie.ea2004}. In other words,
despite the variational equations being nonlinear, they can still be
solved exactly! In fact, the variational equations (\ref{kkt}) are
piece-wise linear and can be solved by starting from $\bar
x(\lambda=\lambda_\mathrm{max})=0$ and letting $\lambda$ decrease:
at each stage where a component of $\bar x(\lambda)$ becomes nonzero
(or, in some exceptional cases, where it becomes zero again) a
`small' \emph{linear} system has to be solved to ensure relations
(\ref{kkt}); this procedure is stopped when the desired value of
$\lambda$ is reached (or when some other stopping criterium such as
e.g. $\|\bar x\|_1=R$ or $|\mathrm{supp}\,(\bar x)|=N$, \ldots\ is
satisfied).

Thus, the computational burden of such an algorithm consists of two
contributions at every node: calculation of the remainder
$K^T(y-Kx)$, and the solution of a linear system of size equal to
the support size of that minimizer. The first contribution is the
same at every step, the second is very small at the start (support
size starts at 0) but dominates after a certain number of steps.

Such an algorithm is implemented in this Mathematica package under
the name \texttt{FindMinimizer}. If $K$ and $y$ are given in terms
of exact numbers (integer, rational, \ldots), this algorithm will
yield an exact answer. If approximate numbers are used as input,
then the algorithm will output a numerical solution.

Other implementations exist (e.g. in Matlab
\cite{Donoho.Stodden.ea2007,IMM2005-03897} or in R
\cite{Efron.Hastie.ea2004}), but they only work with approximate,
floating-point, numbers. Being able to work with exact arithmetic is
advantageous for pedagogical reasons: One can easily follow how the
solution is found (the implementation can also print out
intermediate results), at least for small matrices.

Another advantage of this implementation is that it is able to
handle an exceptional case that is left unaddressed in
\cite{Efron.Hastie.ea2004,Donoho.Stodden.ea2007,IMM2005-03897}.
Indeed, in \cite{Efron.Hastie.ea2004} it is mentioned that
their implementation of the algorithm cannot handle the case
where two new indices may enter the support of $\bar x$ at the
same time. Using floating-point data, this is indeed a rare
occurrence, but using small integer matrices and data, it is
easy to find such exceptions (see several fully worked-out
examples in section \ref{examplesection}). Section
\ref{examplesection} describes how the present implementation
deals with that case.

Probably the biggest advantage over existing software is that
L1Packv2 can handle the problem of minimization of the functional
\begin{equation}
\bar x(\lambda)=\arg\min_x \|Kx-y\|^2+2\lambda \sum_i w_i|x_i|
\label{weightedl1functional}
\end{equation}
with weights $w_i\geq 0$. If all weights are nonzero, then this
can be reduced to the problem (\ref{xbar}) by a simple
rescaling of the independent variables. However, if zero
weights are present, this is no longer true. Still, the
implementation provided by L1Packv2 handles that case
(\cite{Efron.Hastie.ea2004,Donoho.Stodden.ea2007,IMM2005-03897}
do not). The variational equations now are:
\begin{equation}
\begin{array}{lclcl}
(K^T(y-K\bar x))_i&= & w_i\, \lambda\; \mathrm{sign}(\bar x_i) & \qquad\mathrm{if} & \bar x_i\neq 0 \\
|(K^T(y-K\bar x))_i|&\leq & w_i\, \lambda & \qquad\mathrm{if} & \bar
x_i= 0
\end{array}\label{weightedkkt}
\end{equation}
In case of the presence of unpenalized components ($w_i=0$ for some
$i$), the biggest difference with the previous case is that the
starting point is no longer the origin (and $\lambda_\mathrm{max}$
is no longer $\max_i|(K^Ty)_i|$). Also, for a given minimizer $\bar
x(\lambda)$, the penalty parameter equals
\begin{equation}
\lambda=\max_{i,w_i\neq 0}|(K^T(y-K\bar x))_i/w_i|
\end{equation}
(compare with expression (\ref{lambdaexpr})). Other than that, the
algorithm (broadly speaking) follows the same lines.

In addition to the LARS/homotopy method used in
\texttt{FindMinimizer}, the package also provides five iterative
algorithms for the minimization problems (\ref{xbar}) and
(\ref{xtilde}):
\begin{enumerate}
\item \texttt{ThresholdedLandweber}
\item \texttt{ProjectedLandweber}
\item \texttt{ProjectedSteepestDescent}
\item \texttt{AdaptiveLandweber}
\item \texttt{AdaptiveSteepestDescent}
\end{enumerate}
After a finite number of steps, these yield an approximation of
the minimizer. Whereas the LARS/homotopy method heavily is
derived from the variational equations (\ref{kkt}), the
`thresholded Landweber' method (see formula (\ref{tlw})) is
inspired by the fixed-point equation (\ref{fixedpoint}). The
formulation (\ref{xtilde}) lays at the basis of the remaining
four algorithms 2--5, see formulas
(\ref{projLW})--(\ref{adaptSD}). They use projection on an
$\ell_1$-ball, which can be calculated fairly easily in
practice \cite{Daubechies.Fornasier.ea2007}. Two of these four
algorithms use a path-following strategy (where the minimizer
$\tilde x(R)$ is sought for a whole range $0\leq R\leq
R_\mathrm{max}$, just as the LARS/homotopy method), whereas the
other two only calculate the minimizer $\tilde x(R)$ for a
specific value of $R$ (just as the `thresholded Landweber'
algorithm calculates $\bar x(\lambda)$ for a specific value of
$\lambda$).

All algorithms have an option to generate a list of intermediate
results at each step (i.e. each node for the exact method, and each
iteration for the iterative methods). One can use the following
variables: iterate or node number, time elapsed since start, $\bar
x$, $y-K\bar x$, $K^T(y-K\bar x)$ (without having to make additional
matrix multiplications!) and any function of these. This is useful
for evaluating the behavior of the algorithms, or e.g. for easily
being able to plot the trade-off curve, $(\|\bar x(\lambda)\|_1$ vs.
$\|K\bar x(\lambda)-y\|^2)$ as a function of $\lambda$. An example
is given in section \ref{examplesection}.

\section{Main algorithm}

The main algorithm is called \texttt{FindMinimizer} and has syntax
\begin{quote}
\texttt{FindMinimizer[K\_, y\_, options\_\_\_]}.
\end{quote}
It calculates the minimizer $\bar x(\lambda)$ of
(\ref{l1functional}) with the `exact' algorithm (i.e. homotopy
method/LARS \cite{Osborne.Presnell.ea2000,Efron.Hastie.ea2004}). It
takes the operator $K$ and the data $y$ as input. The option
\texttt{Weights $\rightarrow w$} can be used to set the weights
$w_i$ as in the formulation (\ref{weightedl1functional}). The
default for this option is all weights equal to $1$: $w_i=1$.

This algorithm starts from $\bar x=0$ (assuming no zero weights) and
lets $\|\bar x\|_1$ grow, and $\lambda$ descend until a stopping
condition is met. There are four different stopping criteria
possible:
\begin{enumerate}
\item The option \texttt{StoppingPenalty$\rightarrow\lambda$} will stop the algorithm
when the penalty parameter reaches this value.  In other words,
\begin{quote}
\texttt{FindMinimizer[$K$,$y$,StoppingPenalty$\rightarrow\lambda$]}
\end{quote}
will calculate $\bar x(\lambda)$ as in (\ref{xbar}). The
default is \texttt{StoppingPenalty$\rightarrow0$.}

\item The option \texttt{MaximumL1Norm->$R$} will stop the algorithm
when the $\ell_1$-norm of the minimizer reaches the value $R$. In
other words, the command
\texttt{FindMinimizer[$K$,$y$,MaximumL1Norm$\rightarrow R$]} will
calculate $\tilde x(R)$ as in (\ref{xtilde}). The default is
\texttt{MaximumL1Norm$\rightarrow\infty$} which means that the
default will never cause the algorithm to stop.

\item The option \texttt{MinimumDiscrepancy->$d$} will make the algorithm stop when $\|K\bar
x-y\|^2=d$. The default is \texttt{MinimumDiscrepancy->0}.

\item The option \texttt{MaximumNonZero$\rightarrow N$} will stop at the first
node of $\bar x$ having $N$ nonzero coefficients. The default is
\texttt{MaximumNonZero$\rightarrow\infty$}. Unpenalized components
($i$ for which $w_i=0$) are always included in this count.
\end{enumerate}

The use of \texttt{StoppingPenalty}, \texttt{MinimumDiscrepancy} or
\texttt{MaximumL1Norm} allows the user to stop the algorithm at
exactly that value of $\lambda$, $\|K\bar x-y\|^2$ or $\|\bar
x(\lambda)\|_1$ respectively: the first node that overshoots this
condition is computed, but linear interpolation with the previous
node is used to arrive at \emph{exactly} the value specified. This
is very useful if you only need to find the minimizer $\bar
x(\lambda)$ for one specific value of $\lambda$, $\|K\bar x-y\|^2$
or $\|\bar x\|_1$. There is no need to compile the list of
intermediate nodes and no need for doing the interpolation manually.

Apart form the four specific stopping criteria above, one can also
use a more general criterion:
\texttt{StoppingCondition$\rightarrow$cond}. E.g.
\texttt{StoppingCondition:>(Time>=3)} will make sure the algorithm
stops at the first \emph{node} that is calculated within 3 seconds.
All the variables listed in table \ref{vartable} can be used to
express a stopping condition. The default is \texttt{False} meaning
that it will not cause the algorithm to stop.

Notice that the default behavior will make the algorithm stop
at the last node, for which $\lambda=0$ and $K^T(y-Kx)=0$. When
using floating point data one should take care to specify a
good stopping condition. Using floating point numbers, the
default stopping condition,
\texttt{StoppingPenalty$\rightarrow0$}, most likely will not
work (because $\max_i|(K^T(y-K\bar x))_i|$ will eventually be
$\approx 10^{-16}$ or so, which is still strictly larger than
$0$); In this case, it is best to specify an explicit stopping
condition of type \texttt{StoppingPenalty$\rightarrow\lambda$}
with $\lambda>0$. One could also use
\texttt{MaximumL1Norm$\rightarrow R$} if one has a good idea
for a suitable final $\|\bar x\|_1$.

The output of the \texttt{FindMinimizer} algorithm is of the form
$\{\{\mathrm{collected \ data} \},\bar x\}$, or just $\bar x$ if no
data are collected (default). The data that are to be collected at
each node are described by the \texttt{ListFunction} option. Table
\ref{vartable} lists the variables that can be used. E.g.
\begin{quote}
\texttt{FindMinimizer[$K$,$y$,\\
\rule{2cm}{0mm}ListFunction:>\{Norm[Minimzer,1],Norm[DataMisfit]${}^2$\}]}
\end{quote}
will make a list containing $\{\|\bar x\|_1,\|K\bar x-y\|^2\}$ at
each node. This list could then be used to plot the trade-off curve
(see section \ref{examplesection} for an example).

\begin{table}
\begin{tabular}{ll}
Name & Meaning \\ \hline
\texttt{Minimizer} & $\bar x(\lambda)$ \\
\texttt{Counter} & index of the node (if $w_i\neq 0$, the zeroth node is $\bar x=0$)\\
\texttt{DataMisfit} & $y-K\bar x$\\
\texttt{Remainder} & $K^T(y-K\bar x)$\\
\texttt{Penalty} & $\lambda$\\
\texttt{Support} & the support of $\bar x$\\
\texttt{Time} & elapsed time since start (in seconds, but no
`\texttt{Seconds}' included)
\end{tabular}
\caption{A list of the variables that can be used in the
\texttt{ListFunction} and \texttt{StoppingCondition} options of the
\texttt{FindMinimizer} function.}\label{vartable}
\end{table}

Finally, there is an option \texttt{Verbose $\rightarrow n$} that
controls the amount of information that is printed while running.
\texttt{Verbose $\rightarrow 0$} corresponds to no information at
all.

The \texttt{FindMinimizer} function temporarily switches off
Mathematica division-by-zero warnings while it is running.

The implementation in SparseLab (called \texttt{SolveLasso})
\cite{Donoho.Stodden.ea2007} can also calculate the minimizers
with only positive components. This is not possible in
\texttt{L1Packv2}.

\section{Other algorithms and utilities}

Firstly, the following auxiliary functions are available.
\begin{itemize}
\item \texttt{SoftThreshold[x\_, $\lambda$\_]}:\\
Soft-thresholding operation of $x$ with threshold $\lambda\geq 0$.
It is defined as:
\begin{equation}
S_\lambda(x)= \left\{\begin{array}{lcl}x-\lambda & &
x>\lambda\\
0 & \qquad & |x|\leq \lambda\\
x+\lambda & & x<-\lambda
\end{array}\right.\label{soft}
\end{equation}
If $x$ is a vector, $S_\lambda$ is applied component-wise. Here
$\lambda$ can also be a vector.

\item \texttt{ProjectionOnL1Ball[x\_, R\_]}:\\
Projection of $x$ on a $\ell_1$ ball of radius $R$. The projection
$P_R(x)$ of $x$ on the $\ell_1$-ball with radius $R$ is defined as:
\begin{equation}
P_R(x)=\arg\min_{\|u\|_1=R} \|x-u\|^2
\end{equation}
One can show that $P_R(x)=S_\tau(x)$ where $\tau$ has to be chosen
as a function of $R$ and $x$ \cite{Daubechies.Fornasier.ea2007}.

\item \texttt{MinimizerInterpolationFunction[minlist\_,lambdalist\_,var\_]}
takes as input a list of nodes of $\bar x$ and the corresponding
value of the parameter $\lambda$ and returns $\bar x(\lambda)$
evaluated at \texttt{var} (can be symbolic). $\bar x(\lambda)$ is a
piece-wise linear function of $\lambda$, so the result of this
procedure is a list of \texttt{InterpolatingFunction}s evaluated at
\texttt{var}. The input range of the function is between the
smallest and the largest $\lambda$ in \texttt{lambdalist}

\item
    \texttt{CheckMinimizerList[K\_,y\_,minlist\_,lambdalist\_]}
    checks whether the points in \texttt{minlist} are
    \emph{consecutive} nodes of the minimizer (\ref{xbar})
    (corresponding to $\lambda$ in \texttt{lambdalist}).
    I.e. it checks the fixed point equation
    (\ref{fixedpoint}) \emph{between} these nodes. It is
    rather slow because it uses \texttt{Simplify} and has
    to recompute the remainders. The result is
    \texttt{True} if the test is successful, \texttt{False}
    if the test fails and \texttt{Indeterminate} if the
    input is not exact (i.e. floating point), not of equal
    length or not sorted in order of descending $\lambda$.
    The only options are \texttt{Weights} and
    \texttt{Verbose}.

\end{itemize}

In addition to the \texttt{FindMinimizer} algorithm that uses the
LARS/homotopy method to calculate the minimizer
(\ref{l1functional}), there are also five \emph{iterative
algorithms} available. They can be used to find an approximation of
the minimizer (\ref{xbar}) or (\ref{xtilde}). The first three
calculate the minimizer for one specific value of the penalty
parameter $\lambda$ in (\ref{xbar}) or one specific value of $R$ in
\ref{xtilde}.
\begin{itemize}
\item \texttt{ThresholdedLandweber[K\_,y\_,lambda\_,options\_\_\_]
} implements the thresholded Landweber algorithm
\cite{Daubechies.Defrise.ea2004}:
\begin{equation}
x^{(n+1)}=S_\lambda[x^{(n)}+K^T(y-Kx^{(n)})]\label{tlw}
\end{equation}
Options are listed in table \ref{iterationtable}. The main advantage
of this algorithm is that it is very simple. Convergence can be
(very) slow, especially for small values of the penalty parameter
$\lambda$.

\item \texttt{ProjectedLandweber[K\_,y\_,R\_options\_\_\_]}
implements the Projected Landweber algorithm
\cite{Daubechies.Fornasier.ea2007}:
\begin{equation}
x^{(n+1)}=P_R[x^{(n)}+\,K^T(y-Kx^{(n)})]\label{projLW}
\end{equation}
Options are listed in table \ref{iterationtable}. This algorithm
uses the $\ell_1$-norm of the minimizer as parameter instead of the
parameter $\lambda$. The nonlinear operator $P_R$ is also slower
than the soft-thresholding operator $S_\lambda$.

\item \texttt{ProjectedSteepestDescent[K\_,y\_,R\_,options\_\_\_]}
implements
\begin{equation}
x^{(n+1)}=P_R[x^{(n)}+\beta^{(n)}\,K^T(y-Kx^{(n)})]\label{projSD}
\end{equation}
with $r^{(n)}=K^T(y-Kx^{(n)})$ and
$\beta^{(n)}=\|r^{(n)}\|^2/\|Kr^{(n)}\|^2$
\cite{Daubechies.Fornasier.ea2007}. Options are listed in table
\ref{iterationtable}. The main advantage of this algorithm is that,
due to the variable step size $\beta^{(n)}$, faster convergence may
be achieved \cite{Daubechies.Fornasier.ea2007}.
\end{itemize}

\begin{table}
\begin{tabular}{lll}
option & default & description\\ \hline
\texttt{StartingPoint} & $\{0,\ldots,0\}$ & specify the zeroth point in the iteration\\
\texttt{Weights} & $\{1,1,\ldots,1\}$ & similar as for \texttt{FindMinimizer} \\
\texttt{ListFunction} & \texttt{Null} &
\parbox[t]{8cm}{similar as for \texttt{FindMinimizer}}\\
\texttt{StoppingCondition} & \texttt{Counter>=1} &
\parbox[t]{8cm}{specify when the algorithm should stop. All
variables that can be used in \texttt{ListFunction} can also be used
here. E.g.: \texttt{StoppingCondition -> (Counter>=50 ||
Time>=10)}}.
\end{tabular}
\caption{Options available for \texttt{ThresholdedLandweber},
\texttt{ProjectedLandweber}, \texttt{ProjectedSteepestDescent}.
Notice that the default behavior will only make one iteration.}
\label{iterationtable}
\end{table}

The last two algorithms try to use an (approximate) path following
strategy by gradually increasing the radius $R$ of the $\ell_1$-ball
on which the minimizer is sought. This can also be interpreted as a
`warm start strategy' in which an (approximate) minimizer (for one
value of $\|\tilde x\|_1$) is used as the starting point for
obtaining the minimizer with a larger value of $\|\tilde x\|_1$.
\begin{itemize}
\item \texttt{AdaptiveLandweber[K\_,y\_,R\_,numsteps\_,options\_\_\_]}
\begin{equation}
x^{(n+1)}=P_{R_n}[x^{(n)}+\,K^T(y-Kx^{(n)})],\qquad
x^{(0)}=0\label{adaptLW}
\end{equation}
with $R_n=(n+1)R/\mathrm{numsteps}$
\cite{Daubechies.Fornasier.ea2007}.

\item \texttt{AdaptiveSteepestDescent[K\_,y\_,R\_,numsteps\_,options\_\_\_]}
\begin{equation}
x^{(n+1)}=P_{R_n}[x^{(n)}+\beta^{(n)}\,K^T(y-Kx^{(n)})],\qquad
x^{(0)}=0\label{adaptSD}
\end{equation}
with $r^{(n)}=K^T(y-Kx^{(n)})$,
$\beta^{(n)}=\|r^{(n)}\|^2/\|Kr^{(n)}\|^2$ and with
$R_n=(n+1)R/\mathrm{numsteps}$ \cite{Daubechies.Fornasier.ea2007}.
\end{itemize}
For \texttt{AdaptiveLandweber} and \texttt{AdaptiveSteepestDescent}
the starting point is always zero.

One may also use linear operators instead of explicit matrices: The
following functions need to be overloaded: \texttt{Dot, Part,
Dimensions, Length}. This is not part of L1Packv2.

The algorithms in this section are not available in SparseLab
\cite{Donoho.Stodden.ea2007}, but the Matlab package also has a
number of algorithms that use a path-following strategy.  It
also includes a number of methods geared towards the problem of
finding $\arg\min_{Kx=y} \|x\|_1$ or even $\arg\min_{Kx=y}
\|x\|_0$.

\section{Examples}

\label{examplesection}

In this section a number of elementary but instructive examples are
presented to clarify the \texttt{FindMinimizer} algorithm, to
illustrate its possibilities, and to make the reader aware of some
special cases. A practical application of sparse regression is also
included, as well as an indication of the problem sizes it is
suitable for.

Roughly speaking (there are exceptional cases), when the weights
$w_i$ all equal 1, the \texttt{FindMinimizer} algorithm performs the
following actions at each step (in accordance with equation
(\ref{kkt})):
\begin{itemize}
\item calculate the remainder $r=K^T(y-K\bar x)$.

\item determine the components of $r$ that are maximal (in absolute value). These will
be the nonzero components in $\bar x$.

\item calculate a walking direction $v$ such that $x+\mu v$ will
lead to a uniform decrease (i.e. equal across the different active
components) of the maximal value of $r$ (this involves solving a
linear system).

\item walk as far as possible (choice of $\mu$) until another component of the
remainder becomes equal in size to the maximal ones, or until one of
the nonzero components of $\bar x$ becomes zero.
\end{itemize}

Generically, only one additional component will enter the support of
$\bar x$ at any given step. In the special case that the set $\arg
\max_i |r_i|\setminus\mathrm{supp}(x)$ has more than one element,
i.e. there is more than one candidate new index, one needs to look
more closely at the walking direction $v$ to determine the right
components to add to the support of $\bar x$. The two conditions
that have to be satisfied are: (1) the sign of the components (in
the support of $\bar x$) of $v$ correspond to the sign of the
components of $r$ (so that $r$ and $\bar x$ will have the same sign
in that component, as prescribed by the equalities in (\ref{kkt}));
(2) the components of $K^TKv$ that are not in the support of $\bar
x$ are at least as large as those who are in the support (so that
the corresponding components of $r$ remain smaller or equal to the
components that are in the support of $\bar x$, as prescribed by the
inequalities in (\ref{kkt}).

Each line in the following tables corresponds to a node in the
minimizer $\bar x(\lambda)$. For each node, some key information is
given: the minimizer $\bar x(\lambda)$, the $\ell_1$-norm $\|\bar
x(\lambda)\|_1$, the penalty $\lambda$, the remainder $K^T(y-K\bar
x)$, the support of $\bar x$ and the size of the support of $\bar
x(\lambda)$. This will allow the reader (in every example) to check
the accuracy of the result only by checking signs and comparing
absolute values of components of the remainder (via equations
(\ref{kkt})).

To save space, vectors are written in row form, and rows of the same
matrix are separated by a semicolon.

The first example is elementary. This can easily be done by hand.
The operator $K$ is the identity (in a five dimensional space) and
the data is: $y=(12, -8, 5, 1, 2)$. The list of nodes of the
minimizer $\bar x(\lambda)$ looks as follows:
\begin{displaymath}
\begin{array}{lllllll}
 \mathrm{Node} & \bar x(\lambda) &
\|\bar x\|_1 & \lambda & K^T(y-K\bar x)  & \mathrm{supp}(\bar x) &
|\mathrm{supp}(\bar x)|\\ \hline
 0 & (0,0,0,0,0) & 0 & 12 & (12,-8,5,1,2) & \{\} & 0 \\
 1 & (4,0,0,0,0) & 4 & 8 & (8,-8,5,1,2) & \{1\} & 1 \\
 2 & (7,-3,0,0,0) & 10 & 5 & (5,-5,5,1,2) & \{1,2\} & 2 \\
 3 & (10,-6,3,0,0) & 19 & 2 & (2,-2,2,1,2) & \{1,2,3\} & 3 \\
 4 & (11,-7,4,0,1) & 23 & 1 & (1,-1,1,1,1) & \{1,2,3,5\} & 4 \\
 5 & (12,-8,5,1,2) & 28 & 0 & (0,0,0,0,0) & \{1,2,3,4,5\} & 5
\end{array}
\end{displaymath}
If one understands this example, one understands the whole
algorithm.

At step zero, the remainder is maximal at component number 1. This
will be the first nonzero component in $\bar x$. The value 12 is the
value of $\lambda$ at this node. Now we walk in the direction of the
first unit vector until another component of the remainder becomes
equal (in absolute value) to the first component. This happens when
$\bar x_1=4$ and it defines the first node. The maximal absolute
value of the components of the remainder is again the penalty
$\lambda$ in this node (i.e. 8). Now we walk in the direction
$(1,-1,0,0,0)$, because it will make the first two components of the
remainder decrease equally (in accordance with equation
(\ref{kkt})). We stop at the point $(7,-3,0,0,0)$ because the third
component of the remainder is then equal to the maximal ones. We
continue in this way until the maximal value of the remainder is
zero.

For $K$ not equal to the identity, linear equations have to be
solved at each step, and this is tedious by hand.

The next example illustrates the fact that when the initial
remainder has two components of equal (maximum) size (first and
third in this case), it does not necessarily follow that the
two corresponding components of $\bar x(\lambda)$ become
nonzero. Here, only the third component of $\bar x$ becomes
nonzero in the first step. We have to wait until the third step
to get a nonzero first component. The matrix is $K=(-3, 4, 4;
-5, 1, 4; 5, 1, -4)$ and the data is\quad $y=(24, 17, -7)$. The
list of nodes is as follows:
\begin{displaymath}
\begin{array}{lllllll}
 \mathrm{Node} & \bar x(\lambda) &
\|\bar x\|_1 & \lambda & K^T(y-K\bar x)  & \mathrm{supp}(\bar x) &
|\mathrm{supp}(\bar x)|\\ \hline
 0 & (0,0,0) & 0 & 192 & (-192,106,192) & \{\} & 0 \\
 1 & \left(0,0,\frac{43}{16}\right) & \frac{43}{16} & 63 &
 \left(-\frac{209}{4},63,63\right) & \{3\} & 1 \\
 2 & \left(0,\frac{43}{15},\frac{43}{15}\right) & \frac{86}{15} &
 \frac{128}{15} & \left(-\frac{128}{15},\frac{128}{15},\frac{128}{15}\right) & \{2,3\} & 2 \\
 3 & \left(-\frac{172}{73},\frac{301}{73},0\right) & \frac{473}{73} &
 \frac{256}{73} & \left(-\frac{256}{73},\frac{256}{73},\frac{256}{73}\right) & \{1,2\} & 2 \\
 4 & \left(-\frac{2356}{991},\frac{4251}{991},0\right) & \frac{6607}{991} &
 \frac{256}{991} & \left(-\frac{256}{991},\frac{256}{991},-\frac{256}{991}\right) & \{1,2\} & 2 \\
 5 & (-4,5,-2) & 11 & 0 & (0,0,0) & \{1,2,3\} & 3
\end{array}
\end{displaymath}
A similar phenomenon may occur anywhere in the calculation (any
step), when the two or more \emph{new} components of the remainder
become equal to the maximum value of the remainder.

In fact, choosing which components are the right ones to enter the
support of $\bar x(\lambda)$ is the hard part of the algorithm.
Indeed, if one knows the support of $\bar x(\lambda)$, one only
needs to solve the linear system in (\ref{kkt}) to find $\bar x$.

When using floating point arithmetic and real data, one is unlikely
to encounter such cases.

Below is the output generated, for the same operator $K$ and data
$y$, by the Matlab implementation \cite{IMM2005-03897}. Clearly, the
first node is wrong, because the largest component of the remainder
(third) does not correspond to the nonzero component in $\bar x$
(first). These points do not satisfy the fixed point equation
(\ref{fixedpoint}); they are not minimizers. The correct minimizers
were given above.
\begin{displaymath}
\begin{array}{lll}
\mathrm{Node} & \bar x & K^T(y-K\bar x)\\ \hline
0&(0,0,0) &(-192.0000,106.0000,192.0000)\\
1&(-1.8298,0,0)&(-84.0426,84.0426,96.8511)\\
2&(-1.0293,0,1.4009)&(-58.4255,71.2340,71.2340)\\
3&(-2.1807,0.6141,0)&(-55.9691,68.7776,68.7776)\\
4&(-2.5971,3.8760,0)&(7.7421,5.0664,-5.0664)\\
5&(-10.3550,7.6758,-9.7765)&(2.6758,-0.0000,0.0000)
\end{array}
\end{displaymath}
These numbers were calculated with the Matlab command
\begin{quote}
\texttt{lars(mat,data,}'\texttt{lasso}'\texttt{,0,0)}
\end{quote}
where \texttt{mat} and \texttt{data} correspond to $K$ and $y$ with
a final all-zero line appended.

The problem stems from the fact that the first remainder has two
equal largest components. The Matlab implementation does not handle
this case. Unfortunately, the Matlab script does not give any
warning.

The SparseLab \cite{Donoho.Stodden.ea2007} implementation
\texttt{SolveLasso} also gives a wrong answer (without warning)
in this case. The first breakpoint returned by
\texttt{SolveLasso} is $x=(5.3750,\,\,\, 0,\,\,\, 9.4062)$; the
corresponding value for $K^T(y-Kx)$ is $(-20,\,\,\, 20,\,\,\,
20)$. The signs of the first component of $x$ and $K^T(y-Kx)$
are unequal, in violation of equation (\ref{kkt}).

\begin{figure}\centering
\resizebox{6.5cm}{!}{\fbox{\includegraphics{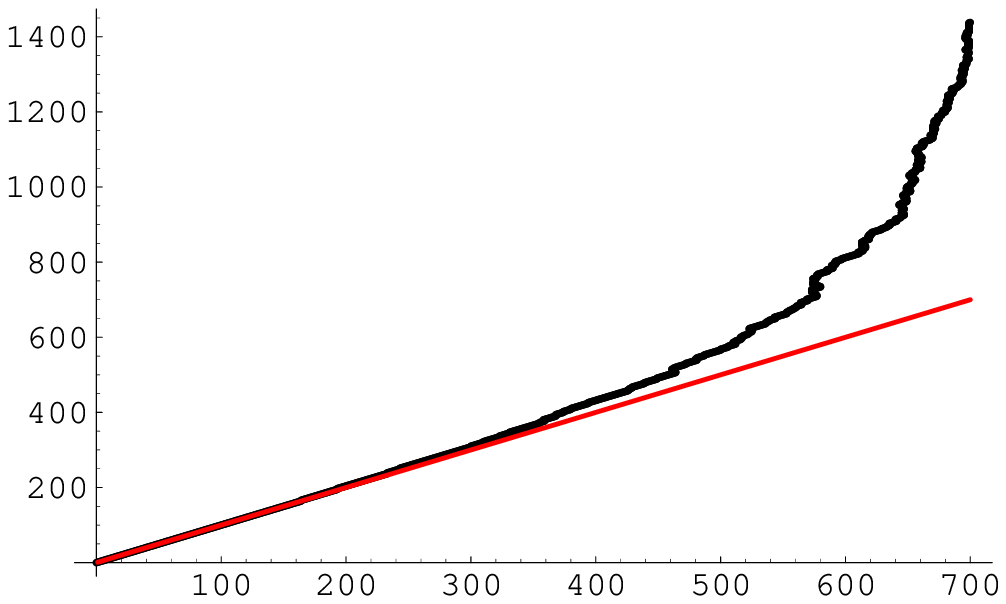}}}
\resizebox{6.5cm}{!}{\fbox{\includegraphics{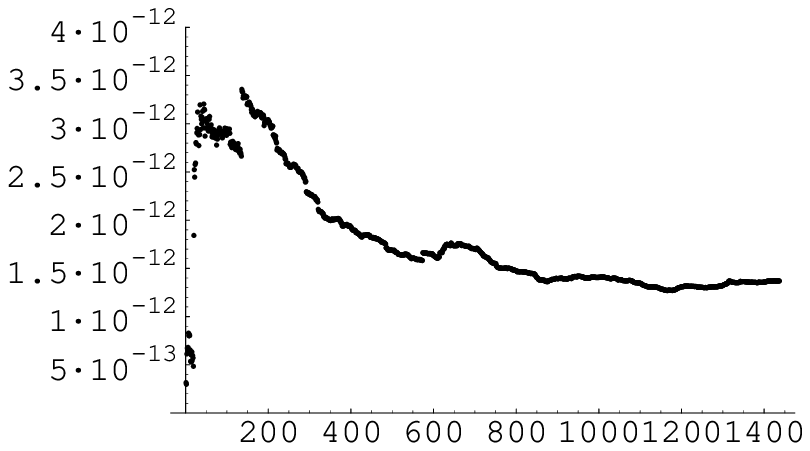}}}
\caption{Left: The node number as a function of the support size of
the minimizer at that node. At first these increase equally
(straight line), but, as the remainder approaches zero, the number
of nodes increases faster than the number of nonzero components in
the minimizer. At the end there are more than twice as many nodes as
there are nonzero components in $\bar x$. Right: The relative error
$\|\bar x-S_\lambda(\bar x+K^T(y-K \bar x))\|/\|\bar x\|$ as a
function of the node number for the same example. Ideally this
should be zero as the minimizers are supposed to satisfy the fixed
point equation (\ref{fixedpoint}). Due to numerical round-off it is
not exactly zero. The matrix $K$ and data $y$ were chosen at random
(see text).}\label{numericfigure}
\end{figure}

The third example shows that a nonzero component of $\bar x$ may be
set to zero again at smaller values of $\lambda$: $K=(-4, 3, -1; -4,
4, 3; -1, 1,-1)$ and $y=(7, 21, 0)$.
\begin{displaymath}
\begin{array}{lllllll}
\mathrm{Node} & \bar x(\lambda) & \|\bar x\|_1 & \lambda &
K^T(y-K\bar x)  & \mathrm{supp}(\bar x) & |\mathrm{supp}(\bar x)|\\
\hline
0 & (0,0,0) & 0 & 112 & (-112,105,56) & \{\} & 0 \\
1 & \left(-\frac{7}{4},0,0\right) & \frac{7}{4} & \frac{217}{4} &
\left(-\frac{217}{4},\frac{217}{4},\frac{175}{4}\right) & \{1\} & 1 \\
2 & \left(0,\frac{7}{3},0\right) & \frac{7}{3} & \frac{133}{3} &
\left(-\frac{133}{3},\frac{133}{3},\frac{112}{3}\right) & \{2\} & 1 \\
3 & \left(0,\frac{49}{18},0\right) & \frac{49}{18} & \frac{308}{9} &
\left(-\frac{595}{18},\frac{308}{9},\frac{308}{9}\right) & \{2\} & 1 \\
4 & \left(0,\frac{28}{9},\frac{7}{3}\right) & \frac{49}{9} & \frac{49}{9} &
\left(-\frac{49}{9},\frac{49}{9},\frac{49}{9}\right) & \{2,3\} & 2 \\
5 & (-1,2,3) & 6 & 0 & (0,0,0) & \{1,2,3\} & 3
\end{array}
\end{displaymath}
It follows that the number of nodes and the size of the support do
not grow at the same rate. In figure \ref{numericfigure} there is a
numerical example that illustrates this phenomenon more clearly.

The final example illustrates the uses of different weights for the
different terms in the penalty as in expression
(\ref{weightedl1functional}). The matrix and data are the same as in
the second example and the weights are $w=(2,1,0)$. The list of
nodes now looks like:
\begin{displaymath}
\begin{array}{lllllll}
\mathrm{Node} & \bar x(\lambda) & \|\bar x\|_1 & \lambda &
K^T(y-K\bar x)  & \mathrm{supp}(\bar x) & |\mathrm{supp}(\bar x)|\\
\hline
 0 & \left(0,0,-\frac{17}{26}\right) & \frac{17}{26} & \frac{214}{13} &
 \left(\frac{659}{26},-\frac{214}{13},0\right) & \{3\} & 1 \\
 1 & \left(0,-\frac{197}{44},\frac{47}{44}\right) & \frac{61}{11} &
 \frac{75}{11} & \left(\frac{150}{11},-\frac{75}{11},0\right) & \{2,3\} & 2 \\
 2 & (3,-4,1) & 8 & 0 & (0,0,0) & \{1,2,3\} & 3
\end{array}
\end{displaymath}
When there is a zero weight (as is here the case for the third
component of $x$), the starting point is no longer necessarily the
origin: the component with zero weight is nonzero. The corresponding
component of the remainder must equal zero at every step. To
determine the penalization parameter $\lambda$ we have to look at
the remainder divided by the weights (component by component
division, and excluding the components with zero weights). In
$\left(\frac{659}{26},-\frac{214}{13},0\right)/(2,1,0)$, the second
component is the largest, and this one will enter the support of
$\bar x$ in the next step.

These tables were, of course, generated automatically with the help
of \texttt{FindMinimizer} function and the built-in Mathematica
\texttt{TeXForm} and \texttt{MatrixForm} commands.

The remaining part of this section is devoted to a small sparse
regression problem.

\begin{figure}\centering
\resizebox{4.5cm}{!}{\fbox{\includegraphics{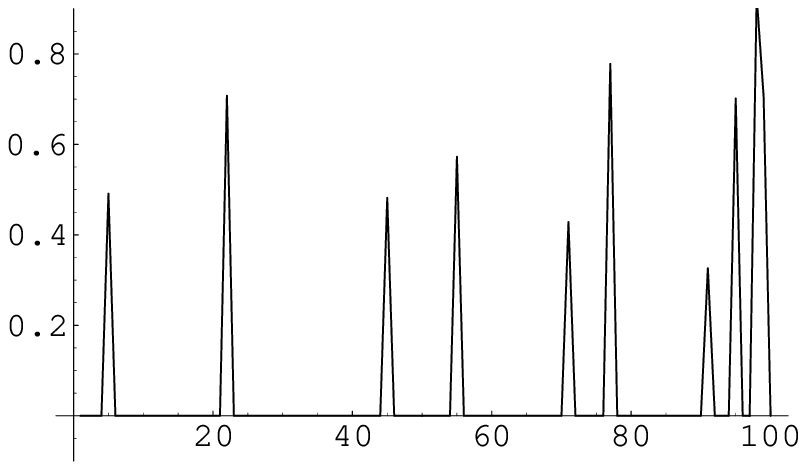}}}
\resizebox{4.5cm}{!}{\fbox{\includegraphics{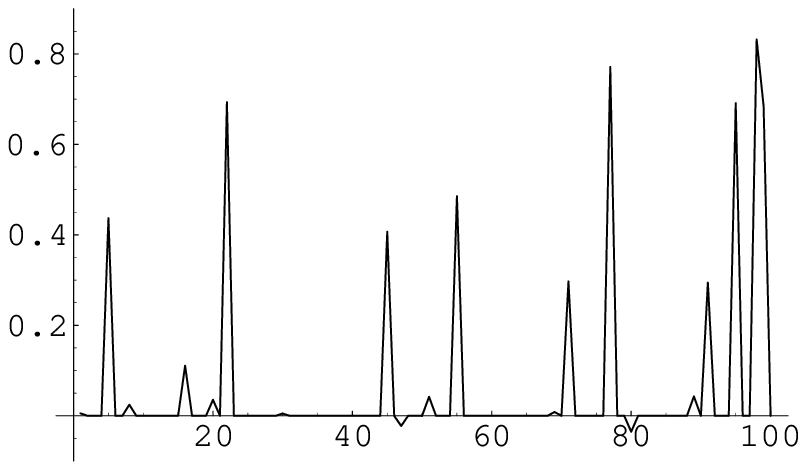}}}
\resizebox{4.5cm}{!}{\fbox{\includegraphics{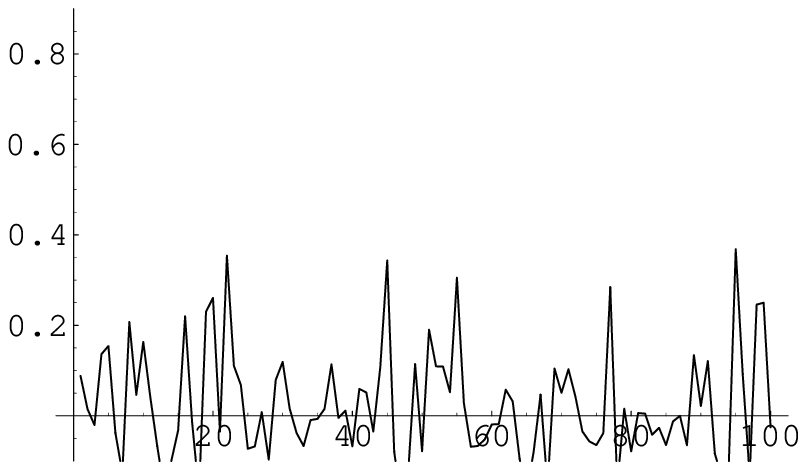}}}
\caption{Left: The (sparse) input model $x_\mathrm{in}$. Middle: the
(sparse) $\ell_1$-penalized reconstruction $x^{(\ell_1)}$. Right:
The $\ell_2$-penalized reconstruction $x^{(\ell_2)}$.
}\label{reconstructionfig}
\end{figure}

We prepare a random matrix $K$ of size $30$ by $100$ with elements
taken from a normal distribution. We also prepare a sparse input
model $x_\mathrm{in}$ with 10 nonzero components (see figure
\ref{reconstructionfig} left). We calculate data $y$ as
$y=Kx_\mathrm{in}+e$ where $e$ is a noise vector with elements taken
from a normal distribution. The norm of the noise vector is 3\% of
the norm of the noiseless data: $\|e\|=0.03 \|Kx_\mathrm{in}\|$. In
other words, we construct an underdetermined toy regression problem
with $100$ unknowns, $30$ noisy data, and $10$ nonzero coefficients
in the sought after solution.

We compare two reconstruction methods: the $\ell_1$ penalized method
(\ref{l1functional}) and the $\ell_2$ penalized method
(\ref{l2functional}). In both cases the penalty parameter is chosen
in such a way that the
 discrepancy of the solution  equals the noise level: $\|K
 x-y\|^2=\|e\|^2$.

For the $\ell_1$-penalized method we can simply use the
\texttt{FindMinimizer} algorithm:
\texttt{FindMinimizer[$K$,$y$,MinimumDiscrepancy$\rightarrow\|e\|^2$]}.
For the $\ell_2$-penalized reconstruction we solve the linear system
(\ref{linsys}), where $\lambda$ is chosen manually in such a way as
to have  $\|K x-y\|^2=\|e\|^2$.

These two reconstructions, $x^{(\ell_1)}$ and $x^{(\ell_2)}$, are
thus equivalent from the point of view of the data: they fit the
data equally well, and no better than the noise level justifies.
However if we compare these two reconstruction with the sparse input
model $x_\mathrm{in}$, we can see that the $\ell_1$ method clearly
succeeds in producing a better reconstruction than the $\ell_2$
method. This toy example thus demonstrates the so-called compressed
sensing idea, whereby a signal can be faithfully reconstructed from
few measurements provided one knows that the signal is sparse
\cite{Donoh2006,CaRoT2006}.

\begin{figure}\centering
\resizebox{6.5cm}{!}{\fbox{\includegraphics{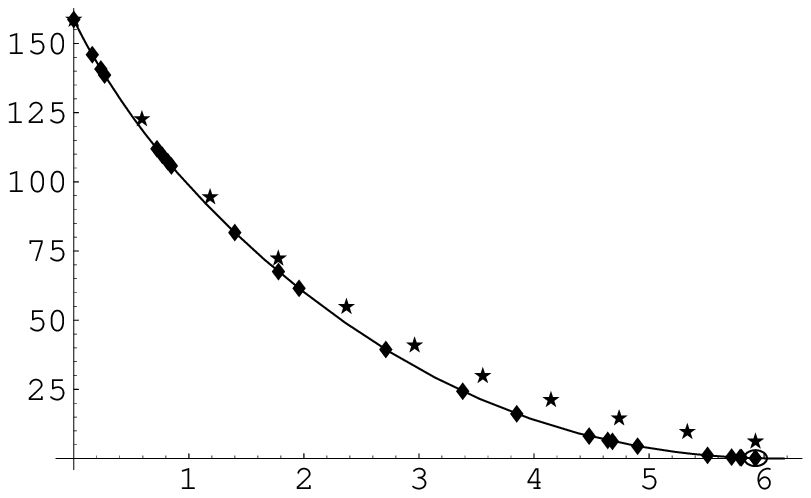}}}
\resizebox{6.5cm}{!}{\fbox{\includegraphics{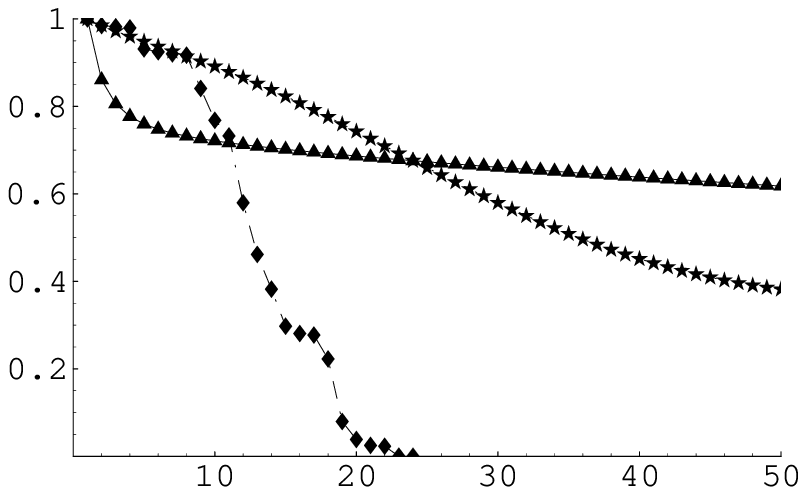}}}
\caption{Left: Trade-off curve (solid line) for the toy regression
problem. Horizontal axis is $\|x\|_1$, vertical axis is
$\|Kx-y\|^2$. Diamonds indicate the nodes of $\bar x(\lambda)$. The
solid line is the trade-off curve (it passes through the diamonds).
The stars represent the path followed by the
\texttt{AdaptiveLandweber} algorithm (\ref{adaptLW}) (10 steps).
Right: Convergence properties of different algorithms. Vertical axis
represents the distance to the true minimizer (i.e. $\|x^{(n)}-\bar
x\|/\|\bar x\|$). Horizontal axis is the iteration number ($n$).
Diamonds represent the nodes of the exact algorithm, triangles are
for the thresholded Landweber algorithm (\ref{tlw}) and stars are
for the adaptive Landweber algorithm \ref{adaptLW}.}\label{pathfig}
\end{figure}

Figure \ref{pathfig} (left) illustrates the trade-off curve for this
problem. The solid line represents the curve $(\|\bar
x(\lambda)\|_1, \|K\bar x(\lambda)-y\|^2)$, for varying $\lambda$.
The points on the trade-off curve can easily be calculated with the
single
command:\\
\texttt{FindMinimizer[$K$,$y$,MinimumDiscrepancy$\rightarrow\|e\|^2$,\\
ListFunction:>\{Norm[Minimizer,1],Norm[DataMisfit]${}^2$\}]}.\\
For
comparison, the approximate trade-off curve traced out by the
adaptive Landweber algorithm (\ref{adaptLW}), 10 steps, is also
pictured. Figure \ref{pathfig} (right) pictures the reconstruction
error $\|\bar x-x^{(n)}\|/\|\bar x\|$ as a function of the iteration
step $n$. Clearly, the iterative algorithms do not perform very
well, and the exact (i.e. LARS/homotopy) method is to be preferred
for this type of problem.

Finally, we make a numerical test of the \texttt{FindMinimizer}
algorithm on a random $700\times 900$ matrix. For large numerical
matrices it is best to always provide an explicit stopping condition
(most likely the user is not interested in the default:
$\lim_{\lambda\rightarrow 0}\bar x(\lambda)$). We keep track of some
intermediate information (node number, time, support size and error)
that is used for plotting in figure \ref{numericfigure}. More than
1400 nodes take just under 8 minutes (on a 2GHz workstation).
\begin{verbatim}
numdata = 700; numvars = 900;
mat=Table[Random[Real,{-3,3}],{numdata}, {numvars}];
data=Table[Random[Real,{-3,3}],{numdata}];
sol=FindMinimizer[mat,data,ListFunction :> {Counter, Time,
Total[Abs[Sign[Minimizer]]],Norm[Minimizer -
SoftThreshold[Minimizer+Remainder, Penalty]]/Norm[Minimizer]}];
\end{verbatim}
These results were plotted in figure \ref{numericfigure}. On the
left we see that the number of nodes and the number of nonzero
components grows equally fast at first. Later on however, as the
algorithm progresses, the support size does not grow with each node
anymore (some nonzero coefficients are set to zero again). On the
right we see that the relative error $\|\bar x-S_\lambda(\bar
x+K^T(y-K \bar x))\|/\|\bar x\|$ (cfr. formula (\ref{fixedpoint}))
stays bounded, even after more than one thousand nodes. This shows
that the \texttt{FindMinimizer} algorithm is very well suited to
handle problems for which the sought after minimizer has fewer than
about 1000 nonzero components.

\section{Summary}

A Mathematica package for the minimization of $\ell_1$ penalized
functionals was described.

The \texttt{FindMinimizer} algorithm finds the exact minimizer of
such a functional, as long as exact input is given. For floating
point data, the same algorithm can be used to obtain a (very good)
numerical approximation of the minimizer.

The \texttt{FindMinimizer} algorithm also handles the case of
different penalties for different components: $\lambda\sum_i w_i
|x_i|$ instead of $\lambda \|x\|_1$, including the case where some
of the $w_i$ are zero. \\
In \cite{FigNo2003}, a wavelet basis is used together with the
iterative soft-thresholding algorithm (\ref{tlw}) to perform
denoising. A wavelet decomposition consists of wavelet coefficients
(encoding local detail) and scaling coefficients (encoding large
scale averages). This is an example where it makes very good sense
to penalize different components differently or not to penalize (at
all) the scaling coefficients. A geophysical example of the use of
such a technique in seismology can be found in \cite{LoNDD2007}.

In many ways, this \texttt{FindMinimizer} algorithm is an
improvement over the Matlab implementation
\cite{IMM2005-03897,Donoho.Stodden.ea2007}; the latter cannot
handle exact numbers (integer, rational), cannot handle zero
weights, and sometimes does not give the correct answer (see
the fourth example in section \ref{examplesection}). The
present algorithm can also compile a list of intermediate
results on the fly. The intermediate value of $x$, $y-Kx$ and
$K^T(y-Kx)$ have to be computed anyway; if one would do this
after finishing the algorithm (based on a list of the nodes),
one has to perform a great number of matrix multiplications all
over again thereby doubling the computing time. Also, it might
be impossible to store all the intermediate minimizer nodes,
rendering the latter strategy impossible in certain cases.
Again, this is not possible with
\cite{Donoho.Stodden.ea2007,IMM2005-03897} (it only outputs a
list of all the intermediate minimizers).
\\
Also with regard to the possible stopping conditions, the
present algorithm has more possibilities: \cite{IMM2005-03897}
can stop at a predetermined number of nonzero components, or at
a predetermined 1-norm of the minimizer. \texttt{SolveLasso} in
SparseLab \cite{Donoho.Stodden.ea2007} can stop at the first
breakpoint with penalty parameter or discrepancy below a
predefined level (though no interpolation with the previous
breakpoint is done to arrive at the exact value) or after a
given number of steps. \texttt{FindMinimizer} has a generic
stopping condition but can also stop at a predetermined value
of the penalization parameter $\lambda$, of $\|\bar x\|_1$ or
at a specific value of the discrepancy $\|K\bar x-y\|^2$
exactly. The latter is very useful in practice when using the
so-called \emph{discrepancy principle} \cite{M.C.1987}: the
minimizer should not fit the data better than the level of the
measurement noise.

The \texttt{FindMinimizer} algorithm was tested on hundreds of
thousands of randomly generated matrices and data (both square and
rectangular). The results were verified and no wrong outcomes were
discovered. This is no guarantee that there aren't any bugs left.
Hence, it is recommended that users would check the results returned
by the algorithms by checking equations (\ref{kkt}) or
(\ref{weightedkkt}).

The \texttt{FindMinimizer} function does not give a warning when the
solution is not unique. If using floating point numbers and a
warning about ill-conditioned matrix is generated, this may mean
that the minimizer is not unique.

The worst case scenario for the \texttt{FindMinimizer} procedure
happens when the first remainder $K^Ty$ has all components of the
same size. In this case, the algorithm has to resort to trying all
possibilities to determine the new index/indices (there are $2^n-1$
nonempty subset of $\{1,2,\ldots,n\}$). Such a special case `never'
happens in practice, but it is easy to manually construct such
examples.

Five other useful algorithms, like the iterative soft-thresholding
algorithm (\ref{tlw}) and algorithms that use projection on an
$\ell_1$-ball (see expressions (\ref{projLW}--\ref{adaptSD})) are
also implemented. They also support the use of (zero and nonzero)
weights. In general, after a finite number of steps, these only
yield an approximative minimizer.

For clarity, only examples with very small size matrices were
given in this manuscript. However, the package was also tested
on the geo-physics application discussed in \cite{LoNDD2007},
and \texttt{FindMinimizer} succeeded in performing that
minimization in less than one minute. The package was also used
to evaluate the performance of competing iterative methods
(other than the ones also discussed here) \cite{Loris2007a}.
The comparison of the relative strength and weaknesses of these
different algorithms lies beyond the scope of this
work.\\
Quite a number of other methods for finding the minimizer
(\ref{l1functional}) are available in SparseLab
\cite{Donoho.Stodden.ea2007}. Some are designed to find the
related $\arg\min_{Kx=y}\|x\|_1$ instead of
(\ref{l1functional}). The $\ell_1$-magic package
\cite{Candes.Romberg2005} also provides a number of primal-dual
interior point methods for solving this and related problems in
Matlab; it does not have an implementation for the
Lasso/homotopy method.

Future work consists of extending the exact algorithm to the case of
linearly constrained minimization:
\begin{equation}
\min_x \|Kx-y\|^2+2\lambda \|x\|_1\qquad\mathrm{with}\quad Ax=a,
\end{equation}
i.e. minimization of the same $\ell_1$ penalized least squares
functional as before, but now under additional linear constraint(s)
described by the linear equation(s) $Ax=a$. Although approximative
methods do exist for this problem, a LARS/homotopy-style (exact)
method for constrained minimization was not described in
\cite{Osborne.Presnell.ea2000,Efron.Hastie.ea2004}. It is possible
to solve this problem exactly and a prototype of such an exact
algorithm, based on the \texttt{FindMinimizer} algorithm, has
already been implemented and used in an application of portfolio
construction (in finance) \cite{Brodie.Daubechies.ea2008}.

\section{Acknowledgements}

Most of the work presented in this article was done while the author
was a `Francqui Foundation intercommunity post-doctoral researcher'
at the Universit\'e Libre de Bruxelles. Part of this work was done
as a post-doctoral research fellow for the F.W.O-Vlaanderen
(Belgium) at the Vrije Universiteit Brussel. Financial help was also
provided by VUB-GOA 62 grant. Stimulating discussions with I.
Daubechies and C. De Mol are gratefully acknowledged. The author
thanks the referees for their suggestions which greatly helped
improve the manuscript.

\bibliography{guide}{}
\bibliographystyle{elsart-num}

\end{document}